\documentclass[reqno]{amsproc}

\usepackage{hyperref}

\hypersetup{
colorlinks=true,
linkcolor=blue,
anchorcolor=blue,
citecolor=blue
}
\usepackage[T1]{fontenc}
\usepackage{lmodern}

\usepackage{lmodern}
\usepackage{amsmath,amssymb,amsthm,mathtools}
\usepackage[letterpaper,margin=1.03in]{geometry}
\usepackage{mathrsfs}
\usepackage{graphicx}
\usepackage{color}

\usepackage[T1]{fontenc}
\usepackage[utf8]{inputenc}
\usepackage{lmodern}
\usepackage{amsmath,amssymb,amsthm,mathtools}
\usepackage{booktabs,array,enumitem}
\usepackage{microtype}
\usepackage{xcolor}
\usepackage{xurl}

\usepackage{bookmark}
\usepackage{fancyhdr}

\newtheorem{theorem}{Theorem}[section]
\newtheorem{lemma}[theorem]{Lemma}
\newtheorem{proposition}[theorem]{Proposition}
\newtheorem{corollary}[theorem]{Corollary}

\newtheorem{remark}[theorem]{Remark}

\numberwithin{equation}{section}

\def\Rm{\operatorname{Rm}}

\newcommand{\nc}{\nabla^{c}}
\newcommand{\nb}{\nabla^{b}}
\newcommand{\Tc}{T^{c}}
\newcommand{\Tb}{T^{b}}

\newcommand{\Gtwo}{\mathrm G_2}
\newcommand{\C}{\mathbb C}
\newcommand{\R}{\mathbb R}
\newcommand{\fa}{\mathfrak a}
\newcommand{\fh}{\mathfrak h}

\newcommand{\cyc}{\mathfrak S}
\newcommand{\dd}{\mathrm d}

\newcommand{\LC}{\nabla^{\mathrm{LC}}}

\DeclareMathOperator{\Hol}{Hol}
\DeclareMathOperator{\ad}{ad}
\DeclareMathOperator{\tr}{tr}
\DeclareMathOperator{\Id}{Id}
\DeclareMathOperator{\diag}{diag}
\DeclareMathOperator{\im}{im}
\DeclareMathOperator{\Rea}{Re}
\DeclareMathOperator{\rad}{rad}

\DeclareMathOperator{\SU}{SU}
\DeclareMathOperator{\Der}{Der}

\DeclareMathOperator{\Lie}{Lie}

\DeclareMathOperator{\Hom}{Hom}

\DeclareMathOperator{\Ad}{Ad}

\newcommand{\fg}{\mathfrak g}
\newcommand{\fk}{\mathfrak k}
\newcommand{\fl}{\mathfrak l}
\newcommand{\fm}{\mathfrak m}
\newcommand{\fr}{\mathfrak r}
\newcommand{\fs}{\mathfrak s}
\newcommand{\fq}{\mathfrak q}

\newcommand{\so}{{\mathfrak{so}}}
\newcommand{\su}{\mathfrak{su}}

\def \C {{\mathbb C}}
\def \R {{\mathbb R}}

\def\Ad {{\rm Ad}}

\def\Rm{{\operatorname R}}

\def\ad{{\operatorname{ad}}}

\numberwithin{equation}{section}

\begin{document}

	\title[Holonomy and connections with parallel torsions]
  {Hermitian connections and holonomy with parallel torsions}

\author{Lei Ni}\thanks{ }
\address{Lei Ni. Zhejiang Normal University, Jinhua, Zhejiang, China, 321004; 
 Department of Mathematics, University of California, San Diego, La Jolla, CA 92093, USA}
\email{lni.math.ucsd@gmail.com; leni@zjnu.edu.cn}



\begin{abstract} 
  
  Let $(M^n,g,J)$ be a connected Hermitian manifold of complex dimension
$n\geq2$, and let $\nabla$ be a Hermitian connection with nonzero parallel
torsion. Assuming that its restricted holonomy acts irreducibly on
$T^{1,0}M$, the author proved that the first Ricci curvature must vanish. In this paper 
we derive a Berger like  complete list of its holonomy representations:
the complexified adjoint representation of a compact simple Lie algebra,
and the standard seven-dimensional complex representation of compact
$\Gtwo$. Curvature of such a connection is proven to be parallel. 

  Secondly, without irreducibility assumption,  we derive
an explicit local normal form for the metric whose Chern connection has parallel trosion,  
which includes holomorphic warpings and nontrivial
holomorphic extensions of hyperk\"ahler manifolds. For complete metrics
with bounded Chern curvature, the characters associated with nontrivial
global restricted-holonomy blocks vanish. Compact metrics admit an unwarped local form, are balanced, and have a
unimodular vertical Lie algebra. On their universal covers, suitable
K\"ahler factors split off and the remaining factor has vanishing first
Chern Ricci form. A compact simply connected elliptic bundle over a K3
surface explains why this local classification does not give a global
product classification. This extends the classification result of Zheng and author on Hermitian 
manifolds whose Chern connection is Ambrose-Singer. In particular, the analysis shows that if the metric 
is not K\"ahler its Chern curvature must satisfies $\Rm_{(\cdot), (\cdot)}X=0$ 
for some parallel nonzero vector $X\in T'M$.

In appendix  we include a by-product: On a compact complex homogenous space $M$ with 
finite fundamental group,  any Hermitian metric  whose Chern connection has parallel torsion,  
must be K\"ahler.
\end{abstract}

\maketitle

\section{Introduction}

In the first part we derive a Berger type list for the irreducible restricted holonomy groups of Hermitian connections 
whose torsion is parallel.  Since honolomy group can be rather arbitary by a result of Hano and Ozeki \cite{HO}, 
the additional assumption on the torsion is necessary.   This seems  to (the author) be a first classification result on 
the holonomy groups
of the  connections with torsions.  
Let $(M^n,g,J)$ be a connected Hermitian manifold: $J$ is integrable and
$g$ is positive definite. 
A 
\emph{Hermitian connection} means a real connection on $TM$ satisfying
\[
\nabla g=0,\qquad \nabla J=0.
\]
We use the conventions
\begin{align}
T(X,Y)&=\nabla_XY-\nabla_YX-[X,Y], \label{eq:torsion}\\
\Rm(X,Y)&=[\nabla_X,\nabla_Y]-\nabla_{[X,Y]}. \label{eq:curvature}
\end{align}
If $\nabla T=0$, the first Bianchi identity is
\begin{equation}\label{eq:bianchi}
\cyc_{X,Y,Z}\Rm(X,Y)Z
=
\cyc_{X,Y,Z}T(T(X,Y),Z).
\end{equation}
The metric and the tensors are extended complex bilinearly when
arguments lie in $TM\otimes_\R\C$.

At a fixed point $p$, set
\[
H=\Hol_p^0(\nabla),\qquad
\fh=\operatorname{Lie}(H),\qquad V=T^{1,0}_pM (\mbox{ also denoted as } T'_pM).
\]
Unless explicitly stated otherwise, \emph{irreducible} means irreducible
on the complex vector space $V$. This differs from irreducibility on
the real vector space $T_pM$.

The first Ricci form refers to the trace of the curvature
\[
\rho^\nabla(X,Y)=\sqrt{-1}\,\tr_\C\!\left(\Rm(X,Y)|_{T^{1,0}M}\right).
\]
Its vanishing is independent of the overall sign convention.

The first-Ricci restriction of the author 
\cite[Theorem~1.2]{Ni-JGA} motivates the further study of irreducible holonomy of Hermitian connections. 
For any Hermitian connection, the vanishing
of its first Ricci form is equivalent to
$H\subseteq\mathrm{SU}(n)$
\cite[Proposition~2.3]{Ni-JGA}. In the paper we prove the following stronger result.

\begin{theorem}[Nonzero parallel torsion holonomy]\label{thm:main}
Let $(M^n,g,J)$ be a connected Hermitian manifold with $n\geq2$, and let
$\nabla$ be a Hermitian connection. Suppose that
\[
T\not\equiv0,\qquad \nabla T=0,
\qquad
H\ \text{acts irreducibly on }T^{1,0}_pM.
\]
Then $\nabla \Rm=0$, 
and, up to unitary conjugacy, the restricted holonomy
representation is one of the following:

\begin{enumerate}

  \item
$\Ad(K)$ acting on $\fk_\C$, where $\fk$ is a compact real simple Lie
algebra, $K$ is a connected compact group with Lie algebra $\fk$, and
$n=\dim_\R\fk$;

\item
compact $\Gtwo$ acting on $\C^7$, the complexification of its standard
real representation on $\R^7$.
\end{enumerate}
Every representation on this list occurs on a non-K\"ahler Hermitian
manifold with a Hermitian connection satisfying the hypotheses.
In every case the representation is of real type, and
\[
H\subseteq\mathrm{SO}(n)\subseteq\mathrm{SU}(n)
\]
after a unitary change of basis.
\end{theorem}

Here $\Ad(K)$ denotes the \emph{effective adjoint group}; its center has
already been divided out. For instance, the first adjoint case is
$\mathrm{SO}(3)\subset\mathrm{SU}(3)$, not the faithful standard
representation of $\mathrm{SU}(2)$.

\begin{theorem}[The Chern and Bismut connections]\label{thm:canonical}
Under the irreducibility assumption above:
\begin{enumerate}
\item If $\nabla=\nc$ is the Chern connection and $\nc\Tc=0$, then
$\Tc=0$. Thus the metric is K\"ahler, and the Chern and Levi-Civita
connections agree. In particular, there are no Hermitian manifold whose Chern connection has
nonzero parallel torsion and complex-irreducible restricted holonomy.
\item If $\nabla=\nb$ is the Bismut connection and
$0\ne\Tb$, $\nb\Tb=0$, then only the adjoint representations
in Theorem~\ref{thm:main}\textup{(i)} occur. Every such representation
occurs.
\end{enumerate}
\end{theorem}

Theorems~\ref{thm:main} and \ref{thm:canonical} classify
\emph{restricted holonomy representations}. If we assume the  completeness and simply-connectedness, Theorem \ref{thm:main} together
 with Ambrose-Singer implies that the manifold is homogenous Hermitian. 
 Without these global assumptions, the
conclusion $\nabla \Rm=0$ is local. In \cite{Agricola} one 
can find many references on interesting interactions between the holonomy and geometry including a Simons type 
holonomy theorem defined for the Lie algebra generated by torsions.

In the second part, let $\nabla=\nc$ be its Chern connection,
with torsion $T^c$ and curvature $\Rm^c$. The  standing 
assumption is
\begin{equation}\label{eq:cpt}
                         \nc T^c=0.
\end{equation}
When the meaning is clear we shall omit the $^c$. 
The manifolds under additional assumption that $\Rm^c$ 
is parallel was studied by Zheng and the author in \cite{NZ-CAS, NZ-classification}. A complete classification was obtained: 
the universal cover splits as the product of a complex Lie group with a Hermitian symmetric space.  
Such a  classification is possible since the holonomy system of such a manifold is same as the Hermitian 
symmetric one and the Ricci flatness can be proven for the non-K\"ahler factor. 
Here we elaborate and substantianly extend the arguments to derive a normal form for 
the Hermitian metrics with parallel Chern torsion. The starting point is the obseravation that the torsion defines 
a Lie algebra structure on the tangent space $T'M$.  
The discussion allows $T^c=0$, so all K\"ahler metrics belong
to the class. Hence the conclusion is a local normal form of the metric, rather than a finite list of
homogeneous spaces. We shall prove that every such metric
is locally of the form
\begin{equation}\label{eq:main-overview}
\begin{aligned}
 \theta&=\theta_A+\sum_i e^{2f_i}z_i\beta_i,\\
 g&=Q(\theta,\overline{\theta})
       +\sum_i e^{f_i+\overline{f_i}}k_i.
\end{aligned}
\end{equation}
Here $A$ is a complex Lie group with a fixed Hermitian inner product $Q$
on its Lie algebra $\mathfrak{a}$; the $z_i$ are parallel sections of the subbundle $F$  (which is also a 
bundle of Lie algebra $\mathfrak{a}$) fixed by the holonomy, associated with 
the holomorphic $(2, 0)$ form $T(\cdot, \cdot)$ defined on the $i$-th irreducible subbundle, valued in $F$, 
 the $k_i$ are K\"ahler metrics; the $f_i$ come from
Lie algebra characters; and $\dd\beta_i$ are parallel holomorphic
symplectic forms (with $\beta_i=g(\cdot, \bar{z}_i)$, after the metric is constructed). 
Precise compatibility conditions and the converse
are given in Theorem~\ref{thm:normal-form}.

The result is obtained by detailed analysis of the Lie algebraic structure on $T'M$ defined by the torsion.
We  remark that 
\begin{enumerate}
\item Parallel Chern torsion does not imply parallel Chern curvature.
      Arbitrary curved K\"ahler metrics can occur among the factors.
\item A splitting of the Chern holonomy representation does not
      automatically give a local metric product. Torsion can obstruct
      integrability of the corresponding distributions.
\item The local normal form does not classify global gluing, monodromy,
      lattices, or bundle periods. 
\end{enumerate}

The well known structure equations (cf. Theorem 2.4 of Ch 3 of \cite{KN-I}) shall be used. 
We also use some basic equations of \cite[\S3, equations (3.1)--(3.5)]{NZ-CAS}.

Recall for the Chern connection $\nc$ the torsion and curvature are defined by 
(\ref{eq:torsion}) and (\ref{eq:curvature}).
Write $V=T^{1,0}M=T'M$. We use the complex bilinear extension of $g$ in
curvature identities, and its associated positive Hermitian form on
$V$ in orthogonality statements. The first Chern Ricci form is
\begin{equation}\label{eq:rho}
                 \rho^{(1)}=\sqrt{-1}\,\tr_{\C}(\Rm^c|_V).
\end{equation}
A different sign convention does not affect its vanishing.

For a Lie-algebra-valued one-form $\theta$, our convention is
\[
 \tfrac12[\theta,\theta](X,Y)=[\theta(X),\theta(Y)].
\]
We use the \emph{right} Maurer--Cartan form, so its structure equation
has a plus sign. In \cite{Ni-JGA}, it was proved for any Hermitian connection 
with parallel torsion if the holonomy is irreducible then $\rho^{(1)}=0$.

For the Bismut connection of a pluriclosed Hermitian metric whose Bismut torsion is parallel, 
a classification was obtained in \cite{BPT}. Our normal form result does not require the 
pluriclosedness assumption. 

The torsion algebra analysis implies that if the metric is complete and simply-connected, and 
there is no K\"ahler part, the manifold is a complex Lie group. Moreover, if there is no nonzero 
global parallel holomorphic vector field, the manifold must be K\"ahler. 

In the appendix we collect detailed list of irreducible holonomies for the Hermitian connections with 
parallel torsions, some needed results on Lie algebras and include a proof of vanishing of holomorphic $1$-forms 
on a compact complex homogeneous space with finite fundamental group, which together with the result in 
Corollary \ref{cor:irreducible} implies that any Hermitian metric with paralell torsion on such a space 
must be K\"ahler. Corollary \ref{cor:irreducible} also implies that for such metric on  any complex manifold the 
 Chern curvature with no kernal must be K\"ahler.

\section{Chern curvature and the holonomy of Chern connection}
\label{sec:blocks}

The defining properties of the Chern connection (\cite{Chern}) are 
\begin{equation}\label{eq:chern-types}
 T^c(V,\overline V)=0,\qquad T^c(V,V)\subset V,\qquad
 \Rm^c(V,V)=\Rm^c(\overline V,\overline V)=0.
\end{equation}
The Bianchi identities (cf. Theorem 5.3 of Ch III of \cite{KN-I}), under $\nabla T=0$, read
\begin{align}
 &\, \cyc_{X,Y,Z}\Rm^c(X,Y)Z
       =\cyc_{X,Y,Z}T^c(T^c(X,Y),Z),\label{eq:first-bianchi}\\
 &\, \cyc_{X,Y,Z}\bigl((\nabla_X\Rm^c)(Y,Z)+\Rm^c(T(X,Y),Z)\bigr)
       =0.\label{eq:second-bianchi}
\end{align}

\begin{lemma}\label{lem:kahler-sym}
Under \eqref{eq:cpt}, the tensor
\[
                    [x,y]_T:=T^c(x,y),\qquad x,y\in T'_pM,
\]
is a complex Lie bracket on each tangent space. It is preserved by
Chern parallel transport. Moreover,
\begin{equation}\label{eq:kahler-sym}
 \Rm^c(x,\overline y)z=\Rm^c(z,\overline y)x
       \qquad (x,y,z\in V),
\end{equation}
and $\Rm^c$ has the usual K\"ahler curvature symmetries.
\end{lemma}

\begin{proof}
Apply \eqref{eq:first-bianchi} to three $(1,0)$ vectors. All curvature
terms vanish, giving Jacobi for $T$. Parallelism of $T$ gives preservation
of the bracket.

Apply the same identity to $x,z,\overline y$. Every iterated-torsion
term vanishes by \eqref{eq:chern-types}, and the remaining curvature
terms give \eqref{eq:kahler-sym}. Metric compatibility and conjugation
then give the corresponding symmetry in the barred slots. In
particular,
\begin{equation}\label{eq:barred-sym}
 g(R(x,\overline y)z,\overline w)
       =g(R(x,\overline w)z,\overline y).
\end{equation}
Together with the curvature type and the metric skew-symmetry, these
are the K\"ahler symmetries.
\end{proof}

\begin{lemma}\label{lem:irreducible-algebra}
Let a connected group $H$ act irreducibly by complex automorphisms on a
nonabelian finite-dimensional complex Lie algebra $\fl$.
Then $\fl$ is simple.
\end{lemma}

\begin{proof}
The solvable radical is characteristic, so it is either $0$ or all of
$\fl$. If $\fl$ is not semi-somple, then its radical is nonzero, hence must be $\fl$. 
The nonzero proper characteristic ideal
$[\fl,\fl]$ contradicts irreducibility for it is nonzero because the
algebra is nonabelian, and proper because the algebra is solvable.
Thus $\fl$ is semisimple. A connected automorphism group preserves
each simple ideal, since its action on the finite set of these ideals
cannot permute them nontrivially. Irreducibility implies that there exists exactly one
simple ideal.
\end{proof}

\begin{proposition}[Parallel Chern torsion implies flatness]
\label{prop:chern-flat}
Suppose $\nc\Tc=0$ and the Lie algebra defined by $\Tc$ at one point
is semisimple. Then $\Rm^c=0$.
\end{proposition}

\begin{proof}
Parallel transport makes the torsion algebras isomorphic. Their Killing
form
\[
B(u,v)=\tr(\ad_u\ad_v),
\qquad
\ad_u(v)=\Tc(u,v),
\]
is therefore a nondegenerate, complex-bilinear, $\nc$-parallel tensor.
Consequently curvature action (as the holonomy action)  is skew with respect to $B$:
\begin{equation}\label{eq:killing-skew}
B(\Rm^c(x,\overline y)z,w)+B(z,\Rm^c(x,\overline y)w)=0.
\end{equation}
For fixed $\overline y$, define
\[
C(x,z,w)=B(\Rm^c(x,\overline y)z,w).
\]
Equation~\eqref{eq:kahler-sym} makes $C$ symmetric in its first
two arguments, and \eqref{eq:killing-skew} makes it skew in its last
two. Hence
\begin{align*}
C(x,z,w)&=C(z,x,w)=-C(z,w,x)\\
        &=-C(w,z,x)=C(w,x,z)=C(x,w,z)=-C(x,z,w).
\end{align*}
Thus $C=0$. Nondegeneracy of $B$ gives $\Rm^c(x,\overline y)=0$ on $V$.
Complex conjugation and the property of Chern curvature \eqref{eq:chern-types} give $\Rm^c=0$.
\end{proof}
The argument of the proof will be used again in Section \ref{sec:algebra} again 
for the case of reducible holonomy actions.

\begin{proof}[Proof of Theorem~\ref{thm:canonical}\textup{(a)}]
Assume $\Tc\ne0$. Lemmas~\ref{lem:kahler-sym}
and \ref{lem:irreducible-algebra} make the torsion algebra simple.
Proposition~\ref{prop:chern-flat} then gives $\Rm^c=0$.
A flat connection has trivial restricted holonomy, which is reducible
on a complex vector space of dimension at least two. This is a
contradiction. In complex dimension one the Chern torsion vanishes trivially.
\end{proof}

If zero torsion is allowed, the question reduces to ordinary K\"ahler
holonomy. In the locally irreducible, Berger's theorem implies that nonsymmetric case the familiar
possibilities are $\mathrm U(n)$, $\mathrm{SU}(n)$, and
$\mathrm{Sp}(m)$ with $n=2m$; locally irreducible Hermitian symmetric
spaces supply the symmetric cases \cite{Berger1955}. The flat complex
line is the one-dimensional trivial-holonomy exception.

Work on the universal cover when global parallel bundles are needed.
At a reference point choose an orthogonal decomposition into complex
restricted-holonomy representations, and extend it by parallel
transport:
\begin{equation}\label{eq:hol-decomp}
                  T'M=F\oplus E_1\oplus\cdots\oplus E_s.
\end{equation}
Here $F$ is the entire fixed subbundle of the \emph{global} restricted
holonomy group, and each $E_i$ is a nontrivial complex-irreducible
summand. All these bundles are holomorphic and Chern parallel. The
connection on $F$ is flat; on the universal cover it has a global
parallel unitary frame.

\begin{proposition}\label{prop:independent}
The holonomy Lie algebra splits into independent factors
\begin{equation}\label{eq:hol-ideals}
               \fh=\fh_1\oplus\cdots\oplus\fh_s,
\end{equation}
where $\fh_i$ acts only on $E_i$. Curvature is supported separately
on these blocks:
\begin{equation}\label{eq:curv-blocks}
\begin{aligned}
 \Rm(F,\cdot)&=R(\overline F,\cdot)=0,\\
 \Rm(E_i,\overline E_j)&=0 &&(i\ne j),\\
 \Rm(E_i,\overline E_i)|_{F\oplus\bigoplus_{j\ne i}E_j}&=0.
\end{aligned}
\end{equation}
\end{proposition}

\begin{proof}
Curvature preserves each parallel summand and annihilates $F$.
For $x\in E_i$, $z\in E_j$, $i\ne j$, the two sides of
\eqref{eq:kahler-sym} belong to distinct summands, so both vanish.
For $x\in F$, the right-hand side is zero, giving $\Rm(F,\overline {T'M})=0$.

If $x\in E_i$ and $y\in E_j$ with $i\ne j$, the endomorphism
$\Rm(x,\overline y)$ sends  $E_i$ to itself. Moreover the above 
implies that $\Rm(x, \overline{y})|_{E_k}=0$ for $k\ne i$.  For $z,w\in E_i$,
\eqref{eq:barred-sym} gives
\[
 g(\Rm(x,\overline y)z,\overline w)
      =g(\Rm(x,\overline w)z,\overline y)=0.
\]
Thus $R(E_i,\overline E_j)=0$. The remaining assertions in
\eqref{eq:curv-blocks} follow in the same way.

The Ambrose--Singer holonomy theorem
\cite[\S2, Theorem~2]{AS-holonomy} says that the holonomy Lie algebra
is spanned by curvature endomorphisms transported along paths.
Parallel transport preserves \eqref{eq:hol-decomp}, and each transported
curvature block $\Rm_{E_i, \overline{E}_i}$ acts only on one $E_i$ nontrivially, which we denote as $\mathfrak{h}_i$.  
Hence the ideals in
\eqref{eq:hol-ideals} act independently.
\end{proof}

Since the Chern connection could have nonzero connection, the splitting of the tangent bundle does not 
immediately yield a de Rham decomposition of the Hermitian metric.

\section{The Bismut connection: the adjoint cases}\label{sec:bismut}

The Bismut connection is the unique Hermitian connection for which
$g(\Tb(X,Y),Z)$ is a three-form \cite{Bismut1989}.
We use the following result of \cite{Ni-JGA}:
if $\nb\Tb=0$, then
\begin{equation}\label{eq:ni-input}
[x,y]_{\fl}:=\Tc(x,y)\quad\text{is a complex Lie bracket on }V,
\qquad \nb\Tc=0.
\end{equation}
The bracket assertion is \cite[Theorem~4.1]{Ni-JGA}; parallelness follows
from the algebraic relation between $\Tb$ and $\Tc$ in
\cite[eq.~(4.3)]{Ni-JGA}, since $\nb$ preserves $g$ and $J$.

\begin{proposition}\label{prop:bismut-adjoint}
If $n\geq2$, $\Tb\ne0$, $\nb\Tb=0$, and the restricted Bismut
holonomy is complex irreducible, then
\[
(H,V)\cong(\Ad(K),\fk_\C)
\]
for a compact real simple Lie algebra $\fk$.
In particular, the real holonomy representation is reducible.
\end{proposition}

\begin{proof}
For a Hermitian metric, $\Tb\ne0$ if and only if $\Tc\ne0$.
By \eqref{eq:ni-input} and Lemma~\ref{lem:irreducible-algebra},
$\fl=(V,\Tc)$ is a complex simple Lie algebra. Since $T^c$ is parallel, the Holonomy acts by
automorphisms of the Lie algebra, and every derivation of a complex semisimple Lie
algebra is inner. Therefore
\[
\fh\subseteq\Der_\C(\fl)=\ad(\fl).
\]
Since $\fl$ is simple, it implies that $\ad$ is injective. Set $\fk=\ad^{-1}(\fh)\subset\fl$, using the injectivity of $\ad$.
Note that $H$ acts on $\fl$ as a Lie algebra automorphism. Since the adjoint action of $H$ 
preserves $\mathfrak{h}$ and $\ad_{ax}=a\ad_x a^{-1}$, for $x\in \fk$, $\ad_x\in \mathfrak{h}$, 
$\ad_{ax}=a\ad_x a^{-1}\in \mathfrak{h}$. Hence $\fk$ is invariant under the action of $H$. 
 The nonzero complex subspace $\fk+\sqrt{-1}\fk$ is $H$-invariant, so
irreducibility gives $\fk+\sqrt{-1}\fk=\fl$.

Also $\fk\cap \sqrt{-1}\fk=0$. Indeed, if $X=\sqrt{-1}Y$ with $X,Y\in\fk$, both
$\ad_X$ and $\ad_Y$ are skew-Hermitian, while
$\ad_X=\sqrt{-1}\ad_Y$ is Hermitian. Thus $\ad_X=0$, and $X=Y=0$.
Consequently
\[
\fl=\fk\oplus \sqrt{-1}\fk
\quad\text{as real vector spaces}.
\]
The algebra $\fk$ is compact because it is isomorphic to the unitary
Lie algebra $\fh$, and its complexification is simple. Hence $\fk$
is compact simple and is a compact real form of $\fl$.
The connected holonomy group is its effective adjoint group.
The real representation splits as the two invariant copies
$\fk$ and $\sqrt{-1}\fk$.
\end{proof}

The relevant facts about inner derivations, compact real forms, and
complex simple Lie algebras belong to the standard structure theory;
see, for example, \cite{Knapp2023}. Explicit realizations are given in
Section~\ref{sec:adjoint-examples} below.

\section{Reduction of the general case}\label{sec:general}

\subsection{The center and the distinction between real and complex irreducibility}
The result below is from \cite{Ni-JGA}. We give a short alternate argument.
\begin{lemma}\label{lem:center}
Under the hypotheses of Theorem~\ref{thm:main}, $\fh$ is nonzero
compact semisimple. In particular, the first Ricci form vanishes.
\end{lemma}

\begin{proof}
A central element of $\fh\subseteq\mathfrak u(n)$ acts, by complex
Schur's lemma, as $cJ$ for some $c\in\R$. If $c\ne0$, its
one-parameter subgroup contains $-\Id$. But torsion is invariant
under parallel transport, and invariance under $-\Id$ would give
\[
-T(x,y)=T(-x,-y)=T(x,y).
\]
This contradicts $T\ne0$. Thus the center of $\fh$ is zero.
A compact Lie algebra is reductive, so $\fh$ is semisimple.
It is nonzero because a trivial connected group cannot act
irreducibly on $\C^n$ when $n\geq2$.

Finally, a representation of a semisimple Lie algebra has trace zero:
the trace vanishes on commutators and $\fh=[\fh,\fh]$.
Thus $\fh\subseteq\su(n)$ and $H\subseteq\mathrm{SU}(n)$.
\end{proof}

The reduction uses a result about \emph{real}
irreducibility. The corrected Cleyton--Swann result, in the formulation
of Moroianu--Schwahn \cite[Theorem~4.4]{MoroianuSchwahn2026} states that a metric connection with nonzero
parallel torsion and real-irreducible infinitesimal holonomy has
totally skew torsion.
(The original reference is
\cite[Theorem~5.14]{CleytonSwann2004}.)

If our real holonomy representation were irreducible, this theorem
would make $\nabla$ the Bismut connection, by uniqueness and
integrability of $J$. Proposition~\ref{prop:bismut-adjoint} would then
make its real representation reducible. Therefore it is reducible.

\begin{lemma}[An orthogonal real form]\label{lem:real-form}
Let $H\subseteq\mathrm U(V,J,g)$ act complex irreducibly and real
reducibly. Then there is an invariant real subspace $E$ such that
\[
V_\R=E\oplus JE, \quad E\mathbin{\perp}JE,
\]
and $E$ is absolutely irreducible. If the representation on $V_\R$
is faithful, so is its restriction to $E$.
\end{lemma}

\begin{proof}
Let $P$ be orthogonal projection onto a nonzero proper real invariant
subspace. Write
\[
P=P_++P_-,
\qquad
P_+=\tfrac12(P-JPJ),\qquad
P_-=\tfrac12(P+JPJ).
\]
Both terms are self-adjoint and commute with $H$; $P_+$ commutes
with $J$ and $P_-$ anticommutes with $J$. Complex Schur's lemma gives
$P_+=a\Id$ with $a\in\R$. Put $A=P_-$. The identity $P^2=P$ becomes
\[
A^2+(2a-1)A+(a^2-a)\Id=0.
\]
Its $J$-anticommuting part is $(2a-1)A=0$. Since $P$ is nontrivial,
$A\ne0$, so $a=\tfrac12$ and $A^2=\tfrac14\Id$.
Thus $\sigma=2A$ is an equivariant self-adjoint involution
anticommuting with $J$. Its orthogonal eigenspaces are $E$ and $JE$.

The natural complex-linear map
$E\otimes_\R\C\longrightarrow(V,J)$, $x\otimes(a+\sqrt{-1}b)\mapsto ax+bJx$,
is an equivariant isomorphism. Thus $E$ is absolutely irreducible.
An element acting trivially on $E$ acts trivially on $JE$ as well,
which proves faithfulness.
\end{proof}

Applied to $T_pM$, the lemma gives
\begin{equation}\label{eq:double}
T_pM=E\oplus JE, \quad E\mathbin{\perp}JE,\qquad
\fh\ \text{acts diagonally by two equivalent faithful real
representations}.
\end{equation}
Since $H$ is connected, its action on the orthonormal real form $E$
lies in $\mathrm{SO}(E)$. This  proves the claimed reduction
to $\mathrm{SO}(n)\subseteq\mathrm{SU}(n)$ in Theorem \ref{thm:main}.

\subsection{Curvature is parallel}

For a metric representation $\fh\subseteq\so(\fm)$, let
$\mathcal K(\fh,\fm)$ denote the space of algebraic curvature tensors:
the maps
\[
A:\Lambda^2\fm\longrightarrow\fh
\]
satisfying the ordinary first Bianchi identity
$\cyc A(x,y)z=0$. As it is well known (cf. Proposition 1.1 of Ch V, \cite{KN-I}) these conditions imply the usual pair-interchange
symmetry of the associated four-tensor. We shall apply the result below to $\nabla \Rm$.

\begin{lemma}\label{lem:curvature-module}
Suppose $\fm=E\oplus F, E\perp F$ and $\fh$ preserves both summands, acting
faithfully on each. Then $\mathcal K(\fh,\fm)=0$.
\end{lemma}

\begin{proof}
For $x,y\in E$ and $z\in F$, the $F$-component of
\[
A(x,y)z+A(y,z)x+A(z,x)y=0
\]
is $A(x,y)z$. Faithfulness on $F$ gives $A(E,E)=0$.
Similarly, $A(F,F)=0$.

For $x\in E$, $y\in F$, and arbitrary $u,v\in\fm$, pair symmetry gives
\[
\langle A(x,y)u,v\rangle
=
\langle A(u,v)x,y\rangle=0,
\]
since $A(u,v)$ preserves $E$. Thus $A(E,F)=0$ as well.
\end{proof}

Differentiate \eqref{eq:bianchi}. Because $\nabla T=0$, its right
side has zero covariant derivative. Thus, for each $u$,
$\nabla_u \Rm$ satisfies the ordinary algebraic first Bianchi identity.
It is metric and takes values in the parallel holonomy-algebra
bundle. By \eqref{eq:double} and Lemma~\ref{lem:curvature-module},
\begin{equation}\label{eq:parallel-curvature}
\nabla \Rm=0.
\end{equation}
This proves that the connection is Ambrose--Singer property without assuming
local homogeneity.

\subsection{The generalized holonomy system and its solvable radical}\label{subsec:transvec}

Parallel torsion and curvature give a {\it generalized symmetric holonomy system:}\footnote{Borrow the term from 
\cite{Sim} where for the Levi-Civita connection, the holonomy system, and symetric holonomy system were introduced.}
\[
\fg=\fh\oplus\fm,\qquad \fm=T_pM,
\]
with bracket
\begin{equation}\label{eq:transvection}
[A,B]=[A,B]_{\fh},\qquad [A,x]=Ax,\qquad
[x,y]=-\Rm(x,y)-T(x,y).
\end{equation}
The first and second Bianchi identities and the holonomy invariance
of $T,\Rm$ give Jacobi identity. This is the Nomizu--Ambrose--Singer
construction; see \cite[p.~62]{Nomizu1954},
\cite[Sections~3.2--3.3]{Kostant1960}, and
\cite{AmbroseSinger1958}. It requires no completeness for this
infinitesimal or local use. Since $\Rm$ is parallel, its values span
$\fh$. The action of $\fh$ on $\fm$ is faithful.

Let $\fr=\rad(\fg)$. By the Levi--Malcev containment result in Appendix~\ref{sec:algebra-appendix}, choose
a real (semi simple) Levi factor $\fs$ containing the semisimple subalgebra $\fh$.
Then
\begin{equation}\label{eq:levi-modules}
\fg=\fs\ltimes\fr,\qquad
\fg/\fh\cong(\fs/\fh)\oplus\fr\cong E\oplus E
\end{equation}
as $\fh$-modules. Note as $\mathfrak{h}$ module $E$ is isomorphic to $JE$. In particular $\fr\cap\fh=0$, and complete
reducibility and the irreducibility of $E$ gives
\[
\fr\cong0,\quad E,\quad\text{or }E\oplus E.
\]

\paragraph{Case 1: $\fr\cong E\oplus E$.}
Equation~\eqref{eq:levi-modules} gives $\fs=\fh$.
Relative to $\fg=\fh\ltimes\fr$, the original invariant complement
$\fm$ is the graph of an equivariant map
\[
\phi:\fr\longrightarrow\fh.
\]
To see this let $P_\fh$ and $P_\fr$ be the two projections to the Levi decomposition components. Consider the map 
$\psi= P_\fr|_\fm: \fm\to \fr$. The map is injective since if $X\in \fm$ with $P_{\fr}(X)=0$, it means 
that $X\in \fh$. But $\fh\cap \fm=\{0\}$. Thus $X=0$. Hence $\psi$ is an isomorphism by  dimension consideration.
Hence $\psi^{-1}(X)\in \fm$ which can be expressed as $\phi(X)+X$ with respect to the Levi-decomposition. 
One can check easily that $\phi$ is $\fh$-invariant.

If $\phi=0$, then $\fm=\fr$ is an ideal, so
$[\fm,\fm]_{\fh}=0$. By \eqref{eq:transvection}, $\Rm=0$, hence
$\fh=0$, a contradiction.

Therefore $\Hom_{\fh}(E,\fh)\ne0$. An irreducible submodule of the
adjoint representation of a compact semisimple algebra is a simple
ideal. If any other simple ideal were present, it would act
trivially on $E$, contrary to faithfulness. Hence $\fh$ is simple
and $E\cong\fh$ is its adjoint representation.

\paragraph{Case 2: $\fr\cong E$.}
Since $\fr$ is solvable and irreducible as an $\fh$-module,
$[\fr,\fr]=0$. The invariant complex structure on $\fg/\fh$ is
integrable, so
\[
\fq=\fh_\C\oplus\fm^{1,0}
\]
is a complex Lie subalgebra of $\fg_\C$, with
\begin{equation}\label{eq:q-intersection}
\fq\cap\overline{\fq}=\fh_\C.
\end{equation}
Here the bar denotes conjugation with respect to the real algebra
$\fg$.

Absolute irreducibility of $E$ implies that
$\fq\cap\fr_\C$ is either $0$ or $\fr_\C$. The latter is impossible:
$\fr_\C$ is conjugation invariant and would then lie in
$\fq\cap\overline{\fq}=\fh_\C$. Thus $\fq\cap\fr_\C=0$.
By dimension, projection $\fq\to\fs_\C$ is an isomorphism.
Consequently $\fq$ is the graph of a linear map
\[
f:\fs_\C\longrightarrow\fr_\C,\qquad f|_{\fh_\C}=0.
\]
Because $\fr$ is abelian, closure of the graph under brackets says
that $f$ is a one-cocycle for the $\fs_\C$-module $\fr_\C$.
Precisely since for $(X, f(X))$ and $(X', (X'))$ in $\fq$, 
$$
[(X, f(X)), (X', f(X'))]=([X, X'], X f(X')-X'f(X))
$$
equals to $([X, X'], f([X, X']))$ by the definition of $f$, we checked the co-cycle condition of $f$. 
Whitehead's lemma in the Appendix \ref{sec:algebra-appendix} gives $f(X)=Xv$ for some $v\in\fr_\C$.
Since $f|_{\fh_\C}=0$, the vector $v$ is $\fh$-fixed.
The module $E_\C$ is nontrivial and irreducible, so $v=0$.
It follows that $\fq=\fs_\C$, which is conjugation invariant.
This contradicts \eqref{eq:q-intersection}. Case 2 cannot occur.

\paragraph{Case 3: $\fr=0$.}
Now $\fg$ is semisimple. Choose a Cartan involution fixing the compact
subalgebra $\fh$. Its compact dual $\fg_u$ contains the same $\fh$
\cite[Chapter~VI, Sections~2--3]{Knapp2023},
and
\begin{equation}\label{eq:compact-dual}
\fg_u/\fh\cong E\oplus E
\end{equation}
as real $\fh$-modules. Indeed, the compact and noncompact parts of
the Cartan decomposition are $\fh$-invariant; multiplication by
$\sqrt{-1}$ on the latter preserves their real representation types.

We now use the compact homogeneous-space classification in precisely
the following form:
an effective compact semisimple pair $(\fg_u,\fh)$ with connected
isotropy and isotropy module consisting of two equivalent
nontrivial irreducible real summands is, up to local isomorphism,
\begin{equation}\label{eq:compact-pairs}
(\so(8),\mathfrak g_2)
\quad\text{or}\quad
(\fk\oplus\fk\oplus\fk,\diag\fk),
\end{equation}
where $\fk$ is compact simple.
For the simple ambient algebra, use the corrected Dickinson--Kerr
classification \cite[Section~2, Example~I.16]{DickinsonKerr2024},
which replaces the relevant list in \cite{DickinsonKerr2008}.
The nonsimple alternative, and the combined statement, can be also founded
in \cite[Section~5.2, p.~20]{PulemotovZiller2023}.

All the hypotheses apply here. Effectivity follows from faithfulness
of the action on $\fm$. The compact semisimple algebra $\fh$ integrates
to a closed connected subgroup of the simply connected compact group
with algebra $\fg_u$; thus passing to compact groups and covers
introduces no closure or local-effectivity problem.

In the first pair of \eqref{eq:compact-pairs}, $E$ is the standard
seven-dimensional $\mathfrak g_2$-module. In the second, it is the
adjoint $\fk$-module. Together with Case 1, these are exactly the
representations in Theorem~\ref{thm:main}. Section~\ref{sec:examples}
proves that both alternatives are realized.

\section{Realizations of the holonomies}\label{sec:examples-1}

\subsection{Adjoint holonomy on complex simple Lie groups}
\label{sec:adjoint-examples}

Let $\fk$ be compact simple, let $\fl=\fk_\C$, and let
$\sigma:\fl\to\fl$ be conjugation with fixed algebra $\fk$.
Write $B$ for the complex Killing form. On a complex Lie group $L$
with Lie algebra $\fl$, use the left-invariant Hermitian metric
whose underlying real inner product is
\begin{equation}\label{eq:cartan-metric}
g(U,V)=-\Rea B(U,\sigma V).
\end{equation}
These Cartan-type metrics are also treated in
\cite[Section~3]{PodestaZheng2025}. We include the connection and
curvature calculation to identify the exact holonomy.

For real left-invariant fields $U=a+\sqrt{-1}b$ and $V=c+\sqrt{-1}d$, where
$a,b,c,d\in\fk$, define
\begin{equation}\label{eq:bismut-formula}
\nabla_UV=[U+\sigma U,V]=2[a,V].
\end{equation}
This specifies a connection in the left-invariant framing.
Every $\ad_a$, $a\in\fk$, is skew-adjoint and complex linear,
so $\nabla g=\nabla J=0$. Its torsion is
\begin{align}
T(U,V)
&=[U,V]+[\sigma U,V]+[U,\sigma V] \notag\\
&=3[a,c]+[b,d]+\sqrt{-1}([a,d]+[b,c]). \label{eq:bismut-torsion}
\end{align}
The identity $\ad_U^*=-\ad_{\sigma U}$ gives
\[
g(T(U,V),W)
=
g([U,V],W)+g([V,W],U)+g([W,U],V),
\]
which is alternating. Thus \eqref{eq:bismut-formula} is the Bismut
connection.

A direct calculation gives
\begin{align}
\Rm (U,V)
&=[\ad_{2a},\ad_{2c}]
  -\ad_{2\Rea[U,V]} \notag\\
&=2\ad_{[a,c]+[b,d]}. \label{eq:bismut-curvature}
\end{align}
The tensors in \eqref{eq:bismut-torsion} and
\eqref{eq:bismut-curvature} have constant components in the
left-invariant framing and are invariant under $\Ad(K)$.
Every connection coefficient in \eqref{eq:bismut-formula} lies in
$\ad(\fk)$, so
\[
\nabla T=\nabla \Rm=0.
\]
Because $[\fk,\fk]=\fk$, curvature values span $\ad(\fk)$.
Parallel transport lies in $\Ad(K)$, and therefore
\[
\Hol^0(\nb)=\Ad(K).
\]
The representation on $\fl=\fk_\C$ is complex irreducible.

Finally, a left-invariant holomorphic frame has constant Hermitian
matrix for \eqref{eq:cartan-metric}. Hence the Chern connection is
flat in that frame and its torsion is the negative of the nonzero
complex Lie bracket. The metric is therefore non-K\"ahler.
This proves the existence assertion for every adjoint entry and
completes Theorem~\ref{thm:canonical}\textup{(b)}.

\subsection{The seven-dimensional \texorpdfstring{$\Gtwo$}{G2} representation}
\label{sec:g2}

We use conventions of \cite{Joyce}. 
Let $W=\R^7$ with its Euclidean metric and positive three-form
\[
\varphi=e^{123}+e^{145}+e^{167}+e^{246}
        -e^{257}-e^{347}-e^{356}.
\]
Define $x\times y$ by
$\langle x\times y,z\rangle=\varphi(x,y,z)$.
For $x\in W$, put
\[
A_x(y)=x\times y,\qquad
B_x=\begin{pmatrix}0&x\\ x^{\mathsf T}&0\end{pmatrix}.
\]
Regard $A_x$ as the upper-left skew block in $\so(7,1)$. Here $x^{\mathsf{T}}$ denotes the tranpose of the $x$. 
The standard $\mathfrak g_2$-decomposition is
\[
\fg=\so(7,1)=\mathfrak g_2\oplus A(W)\oplus B(W),
\qquad
\fm=A(W)\oplus B(W).
\]

With the inner product $-\tr$ on $\so(7)$, the cross-product
identities give
\begin{align*}
-\tr(A_xA_y)&=6\langle x,y\rangle,\\
-\tr([A_x,A_y]A_z)&=-6\varphi(x,y,z),\\
-\tr\bigl((xy^{\mathsf T}-yx^{\mathsf T})A_z\bigr)
  &=-2\varphi(x,y,z).
\end{align*}
Thus the projected brackets are
\begin{equation}\label{eq:g2-brackets}
[A_x,A_y]_{\fm}=-A_{x\times y},\qquad
[B_x,B_y]_{\fm}=-\tfrac13A_{x\times y},\qquad
[A_x,B_y]=B_{x\times y}.
\end{equation}
These formulas can also be checked directly from the displayed
$\varphi$ and the block matrices.

Define an invariant almost complex structure and metric by
\begin{equation}\label{eq:g2-j}
JA_x=-3B_x,\qquad JB_x=\tfrac13A_x,
\end{equation}
and
\begin{equation}\label{eq:g2-metric}
g(A_x,A_y)=9\langle x,y\rangle,\qquad
g(B_x,B_y)=\langle x,y\rangle,\qquad
g(A_x,B_y)=0.
\end{equation}
The metric is positive definite and $J$-Hermitian. The vectors
\[
Z_x=A_x+3\sqrt{-1}B_x
\]
have type $(1,0)$, and \eqref{eq:g2-brackets} gives
\[
[Z_x,Z_y]_{\fm}=2Z_{x\times y}.
\]
Therefore $\mathfrak g_{2,\C}\oplus\fm^{1,0}$ is a subalgebra of
$\fg_\C$, proving integrability of $J$.

On the homogeneous space
\[
M=\operatorname{Spin}_0(7,1)/\Gtwo
\]
take the canonical connection of the displayed reductive decomposition.
It preserves the invariant metric and complex structure, and
\begin{equation}\label{eq:canonical-homogeneous}
T(X,Y)=-[X,Y]_{\fm},\qquad
\Rm(X,Y)=-\ad([X,Y]_{\mathfrak g_2})|_{\fm},\qquad
\nabla T=\nabla \Rm=0.
\end{equation}
These are the standard canonical-connection formulas
\cite{Nomizu1954}. The torsion is nonzero by
\eqref{eq:g2-brackets}.

To identify the entire restricted holonomy, observe that
\[
[B_x,B_y]=xy^{\mathsf T}-yx^{\mathsf T}
\]
spans $\so(7)$. Its projections onto $\mathfrak g_2$ span
$\mathfrak g_2$. Hence the curvature endomorphisms in
\eqref{eq:canonical-homogeneous} span the full isotropy algebra,
and
\[
\Hol^0(\nabla)=\Gtwo
\]
in its standard complexified seven-dimensional representation.

This metric is non-K\"ahler. For
$\omega(X,Y)=g(JX,Y)$, the invariant-form formula gives
\[
d\omega(X,Y,Z)
=-\cyc_{X,Y,Z}\omega([X,Y]_{\fm},Z).
\]
Using \eqref{eq:g2-brackets}--\eqref{eq:g2-metric},
\begin{equation}\label{eq:g2-nonkahler}
d\omega(A_x,A_y,B_z)=-9\varphi(x,y,z)\ne0.
\end{equation}
The connection is not Bismut: the two cyclic torsion contractions are
\[
g(T(A_x,B_y),B_z)=-\varphi(x,y,z),\qquad
g(T(B_y,B_z),A_x)=3\varphi(x,y,z),
\]
so its torsion is not a three-form.

For a complementary description of the underlying complex homogeneous
geometry, the proper transitive action of $\operatorname{Spin}(7,1)$
on the complex quadric
$\mathrm{SO}(8,\C)/\mathrm{SO}(7,\C)$ appears in
\cite[Section~1.2, Table~1.1(vii)]
{KasselKobayashi2025}.

The explicit calculation above supplies the metric, the connection,
and the exact holonomy needed here. This completes the realization
part of Theorem~\ref{thm:main}.

\section{The algebra of parallel Chern torsion}
\label{sec:algebra}

From now on $\nabla$ is the Chern connection, $T$ and $\Rm$ are Chern torsion and curvature. We study the 
Hermitian metrics whose Chern connection has parallel torsion without assuming the holonomy action being irreducible.
 Since $T$ is parallel if $X, Y$ are two parallel global sections of $F$ (after lift to the universal covering)
 we have that $T_{X, Y}$ is also parallel. Hence we have $h(T_{X, Y})=T_{X,Y}$ for any element 
 $h$ of the restricted holonomy. Hence we have that $T_{X, Y}$ is also an section of $F$. At 
 the Lie algebra level for any $A\in \mathfrak{h}$, $T$ is parallel implies that $AT(x,y)=T(Ax, y)+T(x, Ay)$. 
 If $x, y\in F$, $A(x)=0=A(y)$, which then implies $AT(x, y)=0$. Hence $T(x, y)\in F$. 
Use a parallel unitary frame to identify $F$ with a fixed Hermitian
vector space $(\fa,Q)$. The restriction of $T$ makes $\fa$ a complex
Lie algebra.

\begin{proposition}\label{prop:torsion-form}  For any $u, v\in \mathfrak{u}$, $x_i, y_i\in E_i$ and $x_j\in E_j$
there are linear functionals $\lambda_i\in\fa^*$ and vectors
$z_i\in\fa$ such that
\begin{align}
 T(u,v)&=[u,v]_{\fa},&
 T(u,x_i)&=\lambda_i(u)x_i,\label{eq:t-vertical}\\
 T(x_i,x_j)&=0\quad(i\ne j),&
 T(x_i,y_i)&=\Psi_i(x_i,y_i)z_i.\label{eq:t-horizontal}
\end{align}
If $z_i\ne0$, then $\Psi_i$ is a parallel nondegenerate complex
alternating form on $E_i$. When $z_i=0$, the last expression is zero
and $\Psi_i$ may be omitted. The linear functionals satisfy 
\begin{align}
 \lambda_i([\fa,\fa])&=0,\label{eq:character}\\
 [u,z_i]_{\fa}&=2\lambda_i(u)z_i
                         &&(u\in\fa),\label{eq:weight}\\
 \lambda_j(z_i)&=0
                         &&\text{for all }i,j.\label{eq:cross-character}
\end{align}
In particular, the span of $\{z_i\}$ is an abelian ideal of $\fa$.
\end{proposition}

\begin{proof} The argument before the statement proves the first equation of (\ref{eq:t-vertical}).
Now for $u\in F$ we define $L_u=T(u, \cdot): E_i\to T'M$. The same argument implies that 
$AL_u(x)=L_u(A(x))$ for any $x\in E_i$. This together with the irreducibility of $E_i$ (with respect to the 
action of $\mathfrak{h}$) and the decomposition of the action in Proposition \ref{prop:independent} 
imply $L_u(E_i)\subset E_i$ and 
Schur's lemma gives the second equation in
\eqref{eq:t-vertical}. The function $\lambda_i(u)$ is a linear functional since 
$$
\lambda_i(u)=\frac{1}{\dim (E_i)}\tr(L_u|_{E_i}).
$$
It is the character function of the adjoint representation of $\mathfrak{a}$ (acting on $E_i$).

The independence in Proposition \ref{prop:independent} gives the mixed vanishing
in \eqref{eq:t-horizontal} by a similar argument since consider $T(x, y)$ with $x\in E_i, y\in E_j$ $i\ne j$. 
Then the arguement shows that for $A\in \mathfrak{h}_i$, $AT(x, y)=T(Ax, y)$. We then have that 
$T(\cdot, y):E_i\to T'M$ is $\mathfrak{h}_i$ invariant. Hence its image must be inside $E_i$. Then we have 
$T(x, y)\in E_i\cap E_j$ which implies $T(E_i, E_j)=0$. 

 For the second equation in (\ref{eq:t-horizontal}), first we claim that 
\[
             T(x_i,y_i)=C_i(x_i,y_i)+B_i(x_i,y_i),
\]
with $C_i:\Lambda^2E_i\to E_i$ and
$B_i:\Lambda^2E_i\to F$ equivariant.
First $T(x_i, y_j)$ has no $E_j$ component $\forall j\ne i$. 
This is due to that if $A\in \mathfrak{h}_j$ we have 
$$
AT(x_i, y_i)=T(A(x_i), y_j)+T(x_i, A(y_j)=0. 
$$
If $\pi_j$ denotes the projection to $E_j$ we have that $A(\pi_j(T(E_i, E_i)))=\pi_j (A(T(E_i, E_i)))=0$. 
Hence $\pi_j(T(E_i, E_i))$ is annilated by $\mathfrak{h}_j$. By the irreducibility we have that it must be zero. 
This proves  $T(E_i, E_i)\subset F\oplus E_i$. 

Any nonzero invariant $2$-form on the irreducible complex
module $E_i$ is nondegenerate since its kernel is invariant. Now it is easy to see that for
 any parallel section $z$ of $F$, $\Psi_z(\cdot, \cdot)=g(T(\cdot, \cdot), \bar{z})$ is a holomorphic $(2,0)$-form such that it equals 
 to  $\partial \xi'$ with $\xi'(\cdot)=g(\cdot, \bar{z})$. The parallelness of $T$ implies that 
 $g(T(\cdot, \cdot), \bar{z})$ is a parallel $(2, 0)$ form. Hence it is invariant under the holonomy 
 over a fixed point.  Then its kernal is an invariant subspace. If $\Psi_z\ne 0$, then the kernal has to be zero. 
 This implies that $\Psi_z$ is nondegenerate $(2, 0)$ form of $E_i$. The invariance of $\Psi_{z'}$ implies that 
 $\Psi_{z'}=c\Psi_z$ by the irreducibility of $E_i$ under the holonomy action (since they can be identified with 
 transformations of $E_i$ which commute with the irreducible actions of $ \mathfrak{h} $).  Hence the space
of such invariant forms has complex dimension at most one, by applying
Schur's lemma to the corresponding maps $E_i\to E_i^*$. Consequently
$B_i=\Psi_i\otimes z_i$ for some parallel $z_i \in F$ when it is nonzero.

The Jacobi identity applied to  $u,x_i,y_i$ gives 
\begin{equation}\label{eq:pre-C}
             \lambda_i(u)C_i=0,\qquad
             [u,z_i]=2\lambda_i(u)z_i.
\end{equation}

If $C_i\ne0$, then $\lambda_i=0$. The $E_i$-component of the Jacobi identity
with three $E_i$ arguments is therefore precisely the Jacobi identity for $C_i$ for the reason that
the $F$-valued terms have no contribution because $T(F,E_i)=0$. Hence $(E_i, C_i)$ is also a Lie algebra.

Suppose $C_i\ne0$. Its radical is characteristic and thus invariant
under the connected irreducible holonomy action. The radical is either
zero or all of $E_i$. The latter would make the nonabelian algebra
solvable, whose nonzero proper derived algebra would contradict
irreducibility. Thus $(E_i,C_i)$ is semisimple. In fact it is simple,
since the connected holonomy group preserves each simple ideal.
Its Killing form $B$ is nondegenerate and parallel.

Fix $\overline y$ for $y\in E_i$. Define a trilinear tensor on $E_i$
\[
                  A(x,z,w)=B(\Rm(x,\overline y)z,w).
\]
It can be checked easily that it is symmetric in $x,z$ by \eqref{eq:kahler-sym}, and skew in $z,w$
because curvature preserves $B$. These symmetries force $A=0$:
\[
\begin{aligned}
 A(x,z,w)&=A(z,x,w)=-A(z,w,x)=-A(w,z,x)\\
         &=A(w,x,z)=A(x,w,z)=-A(x,z,w).
\end{aligned}
\]
Nondegeneracy of $B$ gives $R_i=0$ everywhere. By
Proposition~\ref{prop:independent} and Ambrose--Singer, this would make
the holonomy action on $E_i$ trivial, a contradiction. Hence $C_i=0$.

The Jacobi identity applied to $u,v,x_i$  gives \eqref{eq:character}. The second
identity in \eqref{eq:pre-C} gives  \eqref{eq:weight}. Setting $u=z_i$
in \eqref{eq:weight} gives $\lambda_i(z_i)=0$ when $z_i\ne0$.
The Jacobi identity applied with two $E_i$ arguments and one $E_j$ argument gives
$\lambda_j(z_i)=0$ for $i\ne j$. The conditions are automatic if
$z_i=0$. Finally, by (\ref{eq:weight}) and $\lambda_j(z_i)=0$ we have 
\[
                    [z_j,z_i]=2\lambda_i(z_j)z_i=0.
\]
Now (\ref{eq:weight}) implies that the span  of $\{z_i\}$ is an  ideal. Hence it is an abelian ideal. 
\end{proof}

\begin{corollary}\label{cor:irreducible}
If $T\ne0$, then $F\ne0$. Thus every non-K\"ahler example has a
nonzero local Chern-parallel holomorphic vector field. If the restricted
Chern holonomy is complex irreducible, then $T=0$.
\end{corollary}

\begin{proof}
If $F=0$, the complete list in Proposition~\ref{prop:torsion-form}
gives $T=0$. In complex dimension at least two, an irreducible
representation cannot have a nonzero fixed subspace. In complex
dimension one, Chern torsion vanishes for dimensional reasons.
\end{proof}

The parallel vector fields are global on the universal cover or if $M$ is simply-connected, otherwise it is only 
locally defined over a simply-connected region.

\section{A local normal form}\label{sec:normal}

\subsection{The construction data}

The local normal form of the Hermitian metric whose Chern connection has parallel torsion involves 
the following data.
\begin{enumerate}
\item A finite-dimensional complex Lie algebra $\fa$ with a positive
      Hermitian inner product $Q$, and its corresponding simply connected complex
      Lie group $A$. Let $\theta_A$ be the right Maurer--Cartan form:
      \begin{equation}\label{eq:MC}
                      \dd\theta_A=\tfrac12[\theta_A,\theta_A].
      \end{equation}
      If $\{X_\alpha\}$ is a base of $\mathfrak{a}$ which are also right invariant vector fields with $\{\omega_\alpha\}$ 
      being their dual, $\theta_A=\sum_{\alpha=1}^{\dim(\mathfrak{a})} \omega_\alpha X_\alpha$.  

\item K\"ahler manifolds $(N_i,k_i)$, $1\le i\le s$, considered locally.
\item Characters $\lambda_i\in\fa^*$ and vectors $z_i\in\fa$
      satisfying \eqref{eq:character}--\eqref{eq:cross-character}.
\item For every $i$ with $z_i\ne0$, a parallel nondegenerate
      holomorphic two-form $\Omega_i$ on $N_i$, and a local holomorphic
      one-form $\beta_i$ with
      \begin{equation}\label{eq:primitive}
                              \dd\beta_i=\partial \beta_i=\Omega_i.
      \end{equation}
      The existence of $\beta_i$ (locally) is  by the Dolbeaut lemma \cite{KM}.
\end{enumerate}
For $z_i\ne0$, the K\"ahler metric $k_i$ is locally hyperk\"ahler due to the existence of nondegenerate 
parallel holomorphic $(2, 0)$-form $\Omega_i$.
No symplectic form is required when $z_i=0$. In this case $N_i$ is just a K\"ahler manifold.

Each $\lambda_i$ integrates on $A$ to a holomorphic additive character
$f_i:A\to\C$, normalized by $f_i(e)=0$, with
\begin{equation}\label{eq:df}
                         \dd f_i=\lambda_i(\theta_A).
\end{equation}
Here 
$\dd \lambda_i(\theta_A)=0$ is due to \eqref{eq:character} and the Maurer-Cartan equation, since 
 if  we write $\theta_A=\sum_{\alpha} \omega_\alpha X_\alpha$
$\lambda_i(\theta_A)=\sum \lambda_i(X_\alpha)\omega_\alpha$,
$$
\dd\lambda_i(\theta_A)= \lambda_i(\dd \theta_A)
=\lambda_i(\sum_{\alpha, \beta}\omega_\alpha\wedge \omega_\beta[X_\alpha, X_\beta])=
=\sum_{\alpha, \beta}\omega_\alpha\wedge \omega_\beta\lambda_i([X_\alpha, X_\beta])=0.
$$
All forms and metrics below are pulled back to
$A\times N_1\times\cdots\times N_s$.

\begin{theorem}[Local normal form]\label{thm:normal-form}
For the above data, set
\begin{equation}\label{eq:theta-model}
          \theta=\theta_A+\sum_{i:z_i\ne0}e^{2f_i}z_i\beta_i
\end{equation}
and
\begin{equation}\label{eq:metric-model}
          g=Q(\theta,\overline\theta)
                    +\sum_{i=1}^s e^{f_i+\overline{f_i}}k_i.
\end{equation}
Then $g$ is a positive Hermitian metric with parallel Chern torsion.

Conversely, every Hermitian metric satisfying $\nc T^c=0$ is locally
holomorphically isometric to a metric of the form
\eqref{eq:theta-model}--\eqref{eq:metric-model}.
No holonomy irreducibility, completeness, or compactness is needed
for this local assertion.
\end{theorem}

\begin{remark}\label{rem:nonminimal}
For the converse, the data can be obtained from the global
restricted-holonomy decomposition on the universal cover.
Nevertheless, the displayed parametrization allows arbitrary K\"ahler
$N_i$, including flat or reducible ones. This avoids asserting that a
globally nontrivial holonomy block remains irreducible or nonflat. Flat factors can  also be absorbed into a larger
vertical Lie algebra, so uniqueness of the data is not claimed.
\end{remark}

\subsection{Verification of the construction}

By \eqref{eq:cross-character},
\begin{equation}\label{eq:lambda-theta}
\lambda_i(\theta)=\lambda_i(\theta_A)=\dd f_i.
\end{equation}
Using \eqref{eq:MC}, \eqref{eq:weight}, and $[z_i,z_j]=0$  direct calculation yields 
\begin{equation}\label{eq:theta-structure}
        \dd\theta=\tfrac12[\theta,\theta]
                      +\sum_i e^{2f_i}z_i\Omega_i,
\end{equation}
where terms with $z_i=0$ are omitted. The equation above can be checked easily by writing out $\theta_A=
\sum \omega_\alpha X_\alpha$. Note that 
$[\theta_A, e^{2f_i}z_i \beta_i]=[e^{2f_i}z_i \beta_i]$.

Choose local unitary coframes $\varphi_i$ for $k_i$, and let $\gamma_i$
be their K\"ahler connection matrices. Thus
\[
                  \dd\varphi_i=-\gamma_i\wedge\varphi_i.
\]
Set
\begin{equation}\label{eq:eta}
                              \eta_i=e^{f_i}\varphi_i.
\end{equation}
Together, a $Q$-orthonormal set of components of $\theta$ and the
$\eta_i$ give a unitary coframe for \eqref{eq:metric-model}. Positivity
follows from the invertible block-triangular change of coframe.
Moreover,
\begin{equation}\label{eq:eta-structure}
          \dd\eta_i=-\gamma_i\wedge\eta_i
                                +\lambda_i(\theta)\wedge\eta_i.
\end{equation}

Consider the unitary connection whose matrix in this coframe is
\begin{equation}\label{eq:connection-matrix}
  \gamma= \diag(0,\gamma_1,\ldots,\gamma_s).
\end{equation}
Equations \eqref{eq:theta-structure} and \eqref{eq:eta-structure}
show that its torsion has only $(2,0)$ and conjugate $(0,2)$ parts. Hence the connection is a Chern connection. 
Its torsion is exactly
\eqref{eq:t-vertical}--\eqref{eq:t-horizontal}.

The Lie algebra structure constants, characters, and $z_i$ are
constant in the vertical parallel frame. If $\Omega_i$ has components
$\Omega_{i,ab}$ in $\varphi_i$, then $e^{2f_i}\Omega_i$ has the
same components in $\eta_i$. Those components are parallel for
$\gamma_i$, because $\Omega_i$ is parallel for $k_i$. Hence every
component of $T$ is parallel for \eqref{eq:connection-matrix}.
This proves the first half of Theorem~\ref{thm:normal-form}.

 Th curvature formulas $\dd \gamma+\gamma\wedge \gamma$ gives that 
\begin{align}
 \Rm&=\diag\bigl(0,\pi_1^*\Rm^{k_1},\ldots,\pi_s^*\Rm^{k_s}\bigr),
       \label{eq:curvature-model}\\
 \rho^{(1)}&=\sum_i\pi_i^*\rho^{k_i}.
       \label{eq:ricci-model}
\end{align}
The first follows from the structure equation in the
adapted frames. The curvature blocks
still depend on the warping factors.

\subsection{Proof of the 2nd half of Theorem \ref{thm:normal-form}}

Start with $\nabla T=0$ and the decomposition
\eqref{eq:hol-decomp}. Work in a sufficiently small neighborhood in
the universal cover. Let $\theta$ be the $\fa$-valued holomorphic
one-form dual to a parallel frame of $F$, annihilating the $E_i$.
Let $\eta_i$ be unitary coframes of $E_i$, with connection matrices
$\gamma_i$. Proposition~\ref{prop:torsion-form} gives
\begin{align}
 \dd\theta-\tfrac12[\theta,\theta]
                         &=\sum_i z_i\Psi_i,\label{eq:original-theta}\\
 \dd\eta_i&=-\gamma_i\wedge\eta_i
                         +\lambda_i(\theta)\wedge\eta_i.
                         \label{eq:original-eta}
\end{align}
Here $\Psi_i$ denotes the parallel two-form on $E_i$ extended to
vanish on all the other summands.

The conditions on the characters imply
\[
 \dd\bigl(\lambda_i(\theta)\bigr)
   =\tfrac12\lambda_i([\theta,\theta])
                +\sum_j\lambda_i(z_j)\Psi_j=0.
\]
Thus there are holomorphic functions $f_i$ with
\begin{equation}\label{eq:original-f}
                          \dd f_i=\lambda_i(\theta).
\end{equation}
$\dd f_i$ annilates all $\{E_j\}$ vectors.

The bundle $F$ is integrable: for parallel vertical fields $U,V$,
\[
                             [U,V]=-T(U,V)\in F.
\]
It defines a regular holomorphic foliation locally, with a local
holomorphic quotient $\pi:M\to B$. By
\eqref{eq:curv-blocks}, curvature vanishes whenever a real vertical
vector is inserted.

Choose unitary frames of each $E_i$ that are parallel along the
simply connected local vertical leaves. In these frames
$\iota_U\gamma_i=0$ for every vertical real field $U$. Since by the equation (\ref{eq:curv-blocks}), 
$\iota_U(\dd\gamma_i+\gamma_i\wedge\gamma_i)=0$, Cartan's formula gives
\[
                            \mathcal L_U\gamma_i=0.
\]
Thus the matrices $\gamma_i$ are basic. Rescale $\eta_i$ to define 
\[
                              \varphi_i=e^{-f_i}\eta_i.
\]
Using \eqref{eq:original-eta} direct calculations shows 
\[
                              \dd\varphi_i=-\gamma_i\wedge\varphi_i.
\]
Since $\iota_U\varphi_i=0$ and $\gamma_i$ is horizontal, Cartan's
formula also gives $\mathcal L_U\varphi_i=0$. Therefore the
$\varphi_i$ descend to coframes on $B$.

Their descended connection is unitary and torsion free. It is the
Levi--Civita connection of a K\"ahler metric on $B$. Its block form
preserves the orthogonal summands, giving a local K\"ahler product
\[
                             B=N_1\times\cdots\times N_s.
\]
For every $z_i\ne0$, the form $e^{-2f_i}\Psi_i$ descends to a parallel
nondegenerate holomorphic two-form $\Omega_i$ on $N_i$. Choose
holomorphic $\beta_i$ satisfying $\dd\beta_i=\Omega_i$, and put
\begin{equation}\label{eq:theta-hat}
                    \widehat\theta
                          =\theta-\sum_i e^{2f_i}z_i\beta_i.
\end{equation}
The weight relations and $[z_i,z_j]=0$ cancel the extra terms in
\eqref{eq:original-theta}, leaving
\[
                    \dd\widehat\theta
                              =\tfrac12[\widehat\theta,\widehat\theta].
\]
The holomorphic Maurer--Cartan integration theorem, equivalently the
holomorphic Frobenius theorem, produces a local map $a:M\to A$ with
$a^*\theta_A=\widehat\theta$. Its differential is an isomorphism on
the vertical tangent spaces, so $(a,\pi)$ is a local biholomorphism.

Finally,
$$
               \dd f_i=\lambda_i(\theta)
                             =\lambda_i(\widehat\theta).
$$
Hence the functions $\{f_i\}$ agree with the functions in \eqref{eq:df}, up to additive
constants. Those constants can be absorbed into the chosen base
metrics and forms. Reconstructing the metric from its unitary
coframes gives \eqref{eq:metric-model}.  This completes the proof.
\hfill$\square$

\section{Examples}
\label{sec:examples}

\subsection{Holomorphic warping over an arbitrary K\"ahler manifold}
\label{subsec:warped}

Let $(N,k)$ be any K\"ahler manifold and let $w$ be the standard
coordinate on $\C$. On $\C\times N$ set
\begin{equation}\label{eq:warped}
                             g=|\dd w|^2+e^{w+\overline w}k.
\end{equation}
This is the normal form with $\fa=\C$, $f(w)=w$, and $z=0$.
For lifted holomorphic vector fields $X,Y$ on $N$, and
$W=\partial_w$, direct calculation gives
\[
 \nabla_WX=X,\qquad \nabla_XW=0,\qquad
 \nabla_XY=D^k_XY,\qquad \nabla_WW=0,
\]
where $D^k$ is the K\"ahler connection of $k$. Hence
\[
                   T(W,X)=X,\qquad T(X,Y)=0,\qquad \nabla T=0.
\]
The curvature block is $\Rm^k$, but as a tensor on the warped metric,
\begin{equation}\label{eq:warped-curvature}
                              \nabla_W\Rm=-\Rm
\end{equation}
on that block. Thus, if $k$ is curved, $\nabla \Rm\ne0$.

This shows that in general we do not have the identity $\Rm(T(\cdot,\cdot),\cdot)=0$
 under $\nabla T=0$ since
\[
                    \Rm(T(W,X),\overline Y)=\Rm(X,\overline Y)
\]
may be nonzero. The  vanishing statement in the CAS
setting uses the additional hypothesis $\nabla \Rm=0$; see
\cite[Lemma~3.1, equation (3.8)]{NZ-CAS}.

\subsection{A holomorphic extension of a hyperk\"ahler metric}
\label{subsec:extension}

Let $(N,k,\Omega)$ be locally hyperk\"ahler, where $\Omega$ is a
parallel holomorphic symplectic form. Choose a local holomorphic
one-form $\beta$ with $\dd\beta=\Omega$, and on $N\times\C_z$ set
\begin{equation}\label{eq:extension}
                         \theta=\dd z+\beta,\qquad
                         g=k+|\theta|^2.
\end{equation}
This is the normal form with $\fa=\C$, $\lambda=0$, and a nonzero
central $z_i$.

In local holomorphic coordinates on $N$, write
\[
                  H_a=\partial_a-\beta_a\partial_z,\qquad
                  Z=\partial_z.
\]
Then
\[
               [H_a,H_b]=-\Omega_{ab}Z,\qquad
               T(H_a,H_b)=\Omega_{ab}Z,\qquad T(Z,\cdot)=0.
\]
In the holomorphic frame $(H_a,Z)$ the metric matrix is
$\diag(k,1)$, so its Chern connection is $\pi^*D^k\oplus\dd$.
Consequently
\[
                  \nabla T=0,\qquad
                  \Rm=\pi^*\Rm^k\oplus0,\qquad \rho^{(1)}=0.
\]
When $k$ is nonflat, this is a non-K\"ahler, non-Chern-flat metric.
The horizontal distribution is not integrable because
$\theta([H_a,H_b])=-\Omega_{ab}$. It therefore illustrates precisely
how a Chern-parallel orthogonal splitting can fail to be a metric
product. The tangent torsion algebra is a complex Heisenberg algebra.

\subsection{The two special cases}

If there are no K\"ahler factors, the construction is a complex Lie
group with a right-invariant Hermitian metric, expressed in a
holomorphic invariant frame. Its Chern connection is flat and its
torsion is parallel. Group inversion converts this to the more usual
left-invariant convention.

If $\fa=0$, the metric is simply a K\"ahler product, and the torsion
vanishes. In particular, the classification necessarily contains the
full freedom of arbitrary K\"ahler metrics.

\section{Completeness, bounded curvature, and compactness}
\label{sec:compact}

\subsection{Vanishing of the characters on nontrivial holonomy blocks}

For this subsection, take the canonical decomposition
\eqref{eq:hol-decomp} from the global restricted holonomy, so all flat
connection summands have been included in $F$.

\begin{proposition}\label{prop:bounded}
Assume $g$ is complete and its Chern curvature has globally bounded
norm. Then every character $\lambda_i$ associated with a nontrivial
block $E_i$ vanishes. Consequently the corresponding $z_i$ lie in
the center of $\fa$.
\end{proposition}

\begin{proof}
Pass to the complete universal cover. A constant vector of $\fa$
defines a global Chern-parallel holomorphic field $Z\in F$.
Using $R(Z,\cdot)=0$ in \eqref{eq:second-bianchi} gives, on the
$i$-th curvature block,
\begin{equation}\label{eq:vertical-derivative}
       \nabla_Z\Rm_i=-\lambda_i(Z)\Rm_i,\qquad
       \nabla_{\overline Z}\Rm_i=-\overline{\lambda_i(Z)}\Rm_i.
\end{equation}
The real vector field $U=Z+\overline Z$ is parallel and has constant
norm. A bounded smooth vector field on a complete Riemannian manifold
has a flow defined for all real times. Along any such trajectory,
\begin{equation}\label{eq:curvature-growth}
 |\Rm_i|^2(t)=e^{-4\Rea(\lambda_i(Z))t}|\Rm_i|^2(0).
\end{equation}
Since $E_i$ is a nontrivial global holonomy block, $R_i$ is nonzero
at some point. Starting the trajectory there, boundedness for both
positive and negative time forces $\Rea(\lambda_i(Z))=0$. Replacing
$Z$ by $\sqrt{-1} Z$ forces its imaginary part to vanish as well.
Thus $\lambda_i=0$; equation \eqref{eq:weight} makes $z_i$ central.
\end{proof}

\begin{corollary}\label{cor:compact-form}
Every compact Hermitian manifold with parallel Chern torsion admits,
locally, an unwarped presentation
\begin{equation}\label{eq:compact-form}
\begin{aligned}
  \theta&=\theta_A+\sum_i z_i\beta_i,\qquad z_i\in Z(\fa),\\
  g&=Q(\theta,\overline\theta)+\sum_i k_i.
\end{aligned}
\end{equation}
The same conclusion holds for complete metrics with bounded Chern
curvature.
\end{corollary}

\begin{remark}\label{rem:flat-characters}
The corollary asserts the \emph{existence} of an unwarped presentation
using the global holonomy decomposition. It does not say that every
nonminimal presentation from Theorem~\ref{thm:normal-form} has zero
characters. A flat K\"ahler factor can occur with a nonzero character
in a Chern-flat presentation; it has been absorbed into $F$
for the argument above.

Similarly, a nontrivial global block may be flat on some open sets.
The proof only needs one point where its curvature is nonzero, and
does not assume nonvanishing at every point.
\end{remark}

\subsection{Balancedness and unimodularity}

\begin{proposition}\label{prop:balanced}
A compact Hermitian manifold with parallel Chern torsion is balanced.
In the canonical unwarped local form \eqref{eq:compact-form}, the
vertical Lie algebra $\fa$ is unimodular.
\end{proposition}

\begin{proof}
For $n\ge2$, let $\omega$ be the fundamental two-form and let $\eta$
be the torsion trace, namely the Gauduchon's $1$-form, with the convention
\[
           \eta(u)=\tr_{\C}\bigl(v\longmapsto T(v,u)\bigr).
\]
The standard trace identities are
\begin{align}
 \partial\omega^{n-1}&=-\eta\wedge\omega^{n-1},\label{eq:trace1}\\
 \partial\overline\partial\omega^{n-1}
   &=(\overline\partial\eta+\eta\wedge\overline\eta)
                                       \wedge\omega^{n-1}.
                                       \label{eq:trace2}
\end{align}
These are the identities used in
\cite[\S3, equations (3.4)--(3.5) and the paragraph preceding
Lemma~3.1]{NZ-CAS}; their compactness argument here only needs
$\nabla T=0$.

Since $\eta$ is parallel, it is holomorphic and
$\overline\partial\eta=0$. Integrating $\sqrt{-1}$ times
\eqref{eq:trace2} and applying Stokes' theorem gives
\[
                  0=\int_M
                        \sqrt{-1}\,\eta\wedge\overline\eta\wedge\omega^{n-1}.
\]
The integrand is a fixed positive multiple of $|\eta|^2\omega^n$,
so $\eta=0$. Equation \eqref{eq:trace1} and its conjugate give
$\dd\omega^{n-1}=0$, which is the balanced condition.

For $u\in\fa$, the trace of \eqref{eq:t-vertical} is
\[
                 \eta(u)=-\tr(\ad_u)
                              -\sum_i(\dim_{\C}E_i)\lambda_i(u).
\]
The characters vanish by Proposition~\ref{prop:bounded}, so
$\tr(\ad_u)=0$ for every $u$. This is unimodularity.
In complex dimension one, torsion is zero and the assertions are
immediate.
\end{proof}

In particular, a compact complex surface with parallel Chern torsion
is K\"ahler, since balancedness in complex dimension two means
$\dd\omega=0$.

\begin{remark}[Nilpotent tangent torsion algebras]\label{rem:nilpotent}
Without compactness, if the whole tangent Lie algebra $(V_p,T_p)$
is nilpotent, then all $\lambda_i$ vanish. Indeed,
$\ad_u$ acts on $E_i$ as $\lambda_i(u)\Id$, whereas every adjoint
endomorphism of a nilpotent Lie algebra is nilpotent. Thus the
unwarped form and centrality of the $z_i$ hold in this case as well.
\end{remark}

\subsection{Global K\"ahler factors and the remaining Ricci curvature}

\begin{theorem}\label{thm:global-splitting}
Let $(M,g,J)$ be complete with parallel Chern torsion and bounded
Chern curvature. Its universal cover admits a holomorphic isometric
product decomposition
\begin{equation}\label{eq:global-splitting}
                          \widetilde M=K\times P,
\end{equation}
where $K$ is K\"ahler and $P$ has parallel Chern torsion and
$\rho_P^{(1)}=0$.

One can take $K$ to contain all nontrivial global holonomy blocks
for which $z_i=0$. The factor $P$ is not asserted to be Chern flat,
and no uniqueness or absence of further flat K\"ahler factors in
$P$ is asserted.
\end{theorem}

\begin{proof}
By Proposition~\ref{prop:bounded}, all $\lambda_i=0$. If $z_i=0$,
then the underlying real distribution
$\mathcal E_i=(E_i\oplus\overline E_i)\cap TM$ satisfies
\[
                    T(\mathcal E_i,TM)=0,\qquad
                    \im T\perp\mathcal E_i.
\]
Let $S=\nabla-\LC$. The contorsion formula is
\begin{equation}\label{eq:contorsion}
 2g(S_XY,Z)=g(T(X,Y),Z)+g(T(Z,X),Y)-g(T(Y,Z),X).
\end{equation}
It gives $S_XY=0$ whenever $Y\in\mathcal E_i$, so
$\mathcal E_i$ is Levi--Civita parallel. It also gives $S_X=0$
when $X\in\mathcal E_i$, so the complex structure is parallel in
the factor directions. The complete simply connected de Rham
theorem~\cite{deRham} produces the global holomorphic isometric K\"ahler
product factors collected into $K$.

Every remaining nonflat block has $z_i\ne0$ and hence a parallel
nondegenerate complex alternating form $\Psi_i$. Its curvature
endomorphisms preserve $\Psi_i$ and have zero complex trace.
Curvature vanishes on $F$. Thus the Chern curvature of the
remaining factor has zero trace, proving $\rho_P^{(1)}=0$.
\end{proof}

For a unitary connection, zero curvature trace is equivalent to its
\emph{restricted} holonomy being contained in $\SU(n)$.
On a nonsimply connected manifold, zero trace curvature need not
make the determinant of the full holonomy trivial. Applied to
\eqref{eq:global-splitting}, the theorem gives this special-unitary
reduction for the residual factor; it does not force the Ricci
form of the arbitrary K\"ahler factor $K$ to vanish.

\section{A compact simply connected example }
\label{sec:global-example}

The construction in Section~\ref{subsec:extension} can be globalized
on a nontrivial holomorphic bundle. The necessary bundle existence
is supplied by Biswas--Dumitrescu~\cite[\S6, Theorem~6.1 and its
proof]{BD}; related torus-bundle constructions can also be found in
Goldstein--Prokushkin~\cite[Theorem~1 and \S2]{GP}.

More precisely, there exist a projective K3 surface $S$, a holomorphic
symplectic form $\Omega$ whose periods generate a lattice
$\Lambda\subset\C$, and a compact simply connected principal elliptic
bundle
\[
                         \pi:P\longrightarrow S
\]
with fiber $\C/\Lambda$, carrying a holomorphic connection one-form
$\theta$ satisfying
\begin{equation}\label{eq:bundle-curv}
                              \dd\theta=\pi^*\Omega.
\end{equation}
This is the bundle and one-form constructed in the proof of
\cite[Theorem~6.1]{BD}; see also its Remark~6.2.

Choose a Ricci-flat K\"ahler metric $k$ on $S$, whose existence follows
from Yau's theorem~\cite{Yau}. The holomorphic symplectic form $\Omega$
is parallel for $k$ by the Bochner principle; see also
\cite[\S1]{BD}. Define
\begin{equation}\label{eq:K3-metric}
                              g=\pi^*k+|\theta|^2.
\end{equation}
On a local holomorphic trivialization,
$\theta=\dd z+\beta$ with $\dd\beta=\Omega$, so the explicit
calculation in Section~\ref{subsec:extension} applies. Therefore
\begin{equation}\label{eq:K3-properties}
           \nabla^cT^c=0,\qquad T^c\ne0,\qquad
           \Rm^c=\pi^*\Rm^k\oplus0,\qquad \rho^{(1)}=0.
\end{equation}

The metric is not Chern flat. Its restricted Chern holonomy is
$\SU(2)$ acting on $\C^2\oplus\C$ by the standard representation
on $\C^2$ and trivially on the last line.

The rank-two nontrivial holonomy distribution is the connection
horizontal distribution. Equation \eqref{eq:bundle-curv} shows that
it is not integrable. Hence this metric is not locally a Hermitian
product of a Chern-flat factor and a K\"ahler factor: in such a
product, the unique nontrivial holonomy block in this dimension would
be tangent to the K\"ahler factor.

Since $P$ is compact and simply connected, passing to the
universal cover does not repair this failure. Thus the local normal
form cannot be improved, under $\nabla T=0$ alone, to a universal-cover
being product of a complex Lie group and K\"ahler manifolds (the result of \cite{NZ-classification}).

\appendix
\section{Dimensions and interpretation}\label{sec:remarks}

\subsection{The adjoint list}

The compact simple types and their adjoint dimensions are as follows
(see the simple Lie algebra tables in \cite[Appendix~C]{Knapp2023}).
Here $n$ is the complex dimension of $M$, not its real dimension.

\begin{center}
\renewcommand{\arraystretch}{1.15}
\begin{tabular}{@{}lll@{}}
\toprule
Type of $\fk_\C$ & Range & $n=\dim_\R\fk$\\
\midrule
$A_r$ & $r\geq1$ & $r(r+2)$\\
$B_r$ & $r\geq2$ & $r(2r+1)$\\
$C_r$ & $r\geq2$ & $r(2r+1)$\\
$D_r$ & $r\geq4$ & $r(2r-1)$\\
$G_2$ & & $14$\\
$F_4$ & & $52$\\
$E_6$ & & $78$\\
$E_7$ & & $133$\\
$E_8$ & & $248$\\
\bottomrule
\end{tabular}
\end{center}

Low-rank isomorphisms are understood; in particular $B_2=C_2$.
The standard $\Gtwo$ entry with $n=7$ is separate from the adjoint
$\Gtwo$ entry with $n=14$.

\subsection{Complex dimension one}

The assumption $n\geq2$ excludes a real exception for arbitrary
Hermitian connections. On the two-dimensional nonabelian real Lie
group, choose a left-invariant orthonormal frame $e_1,e_2$ with
$[e_1,e_2]=e_2$, put $Je_1=e_2$, and declare the frame parallel.
The resulting Hermitian connection satisfies
\[
T(e_1,e_2)=-e_2,\qquad \nabla T=0,\qquad \Rm=0.
\]
Every almost complex structure in real dimension two is integrable.
The trivial restricted holonomy group acts irreducibly on a
one-dimensional complex vector space. Conversely, the scalar-center
argument in Lemma~\ref{lem:center} forces this trivial group whenever
$n=1$ and $T\ne0$ is parallel. Thus the additional dimension-one
entry is exactly $H=\{1\}$.

There is no corresponding nonzero-torsion Chern or Bismut example
in complex dimension one.

\subsection{What the classification adds to the first-Ricci conclusion}

Under the nonzero-torsion hypotheses, the chain of conclusions is
\[
\begin{gathered}
\text{complex-irreducible restricted holonomy}\\
\Longrightarrow\quad
\text{real type and }\nabla \Rm=0\\
\Longrightarrow\quad
\text{compact simple adjoint representation or standard }\Gtwo.
\end{gathered}
\]
The first-Ricci vanishing follows from the center argument \cite{Ni-JGA};
the classification additionally determines the representation and
forces local homogeneity of the connection. The essential
Berger-style feature is that the algebraic curvature constraints,
combined with parallel torsion and integrability, reduce 
a representation-theoretic problem to a short list.

\section{The algebraic facts used in the radical argument}
\label{sec:algebra-appendix}

Here it includes   two standard results used in Section \ref{subsec:transvec}.
All Lie algebras and modules in this appendix are finite dimensional,
over $\R$ or $\C$.
The first result is part (i) of Lemma 3.6 of Ch III, part I of \cite{IT}.
\begin{lemma}[Whitehead]\label{lem:whitehead} Semi-simpleness implies that $H^1(\mathfrak{g}, W)=0$. Namely, 
Let $\fs$ be semisimple and $W$ an $\fs$-module. If a linear map
$f:\fs\to W$ satisfies
\[
f([X,Y])=Xf(Y)-Yf(X),
\]
then $f(X)=Xv$ for some $v\in W$.
\end{lemma}

\begin{proof} The above cited reference provides a proof. Below is a different one. 
Give $W\oplus\mathbb F$, where $\mathbb F$ is the base field,  define the action
\[
X\cdot(w,t)=(Xw+t f(X),0).
\]
The cocycle identity makes this a representation of $\mathfrak{g}$ in $\mathfrak{gl}(W\oplus \mathbb{F})$. Complete reducibility
\cite[Theorem~5.29]{Knapp2023} splits off the trivial quotient, so
there is a vector $(w_0,1)$ which annilated by all $X$. Thus $f(X)=-Xw_0$; take $v=-w_0$.
For a real algebra, the same conclusion follows by complexifying and
taking real parts.
\end{proof}

\begin{lemma}[Levi containment]\label{lem:levi-containment}
If $\fh$ is a semisimple subalgebra of a Lie algebra $\fg$, some (semi-simple)
Levi factor of $\fg$ contains $\fh$.
\end{lemma}

\begin{proof}
Use induction on the dimension of $\rad(\fg)$. There is nothing to
prove when the radical is zero. Otherwise choose a nonzero abelian
ideal $\mathfrak a\subseteq\rad(\fg)$, for example the last nonzero
term of its derived series. The intersection
$\fh\cap\mathfrak a$ is zero, since it is a solvable ideal of $\fh$.
By induction, a semi-simple Levi factor of $\fg/\mathfrak a$ contains the image
of $\fh$. Its inverse image $\mathfrak p$ contains $\fh$ and has
radical $\mathfrak a$.

The existence of a Levi decomposition
\cite[Theorem~B.2]{Knapp2023} gives
$\mathfrak p=\fl\ltimes\mathfrak a$.
Projection to $\fl$ identifies $\fh$ with a semisimple subalgebra
$\fh_0\subseteq\fl$, and $\fh$ is the graph of a cocycle
$f:\fh_0\to\mathfrak a$. Lemma~\ref{lem:whitehead} gives
$f(X)=Xv$ for some $v\in\mathfrak a$.
Because $\mathfrak a$ is abelian,
\[
\exp(\ad v)\bigl(X+f(X)\bigr)=X+f(X)-Xv=X.
\]
Thus a conjugate of $\fl$ contains $\fh$. It is a (semi-simple) Levi factor of
$\fg$, since its dimension is $\dim(\fg/\rad(\fg))$.
\end{proof}

The result seems a little different from the conjugacy theorem of Malcev-Harish-Chandra. 

\section{Chern connection with parallel torsion on C-spaces}
\label{sec:C-spaces-appendix}

Let $(M,J,g)$ be a connected Hermitian manifold, and let $\nc$ denote
its Chern connection.  Recall the torsion and curvature are defined by (\ref{eq:torsion}) and (\ref{eq:curvature})
When no ambiguity is possible, we write $\nabla=\nc$, $T=T^c$, and
$\Rm=\Rm^c$. The real metric is extended complex-bilinearly, so the
Hermitian inner product on $T^{1,0}M$ is $h(u,v)=g(u,\bar v)$.

\begin{theorem}\label{thm:main-appendix}
Let $(M,J)$ be a compact simply connected (or with finite fundamental group) complex  homogenous manifold.
For every Hermitian metric $g$ on $(M,J)$,
\begin{equation}\label{eq:main-implication}
                  \nc\Tc=0\quad\Longrightarrow\quad \Tc=0.
\end{equation}
Consequently, $g$ is K\"ahler.
\end{theorem}

 By results of Montgomery and Wang \cite{Wang1954} for any compact complex homogenous manifold $M$ if the fundamental group is finite, 
 one can find a semi-simple compact group $K$ of biholomorphisms acting transitively on $M$. 

The compact homogeneous complex setting includes the compact simply
connected $C$-spaces of Wang; see \cite{Wang1954,NiWallach2026} for the backgrounds. 
Theorem~\ref{thm:main-appendix} does not require the
metric $g$ to be invariant under the compact transitive group.

\begin{corollary}\label{cor:homogeneous}
Let $M=G/H$ be a compact simply connected $C$-space and let $g$ be a
$G$-invariant Hermitian metric. If its Chern connection has parallel
torsion, then $g$ is K\"ahler.
\end{corollary}

The proof is by contradiction of the following two facts:
\begin{equation}\label{eq:proof-ingredients}
\begin{aligned}
&\text{compact simply connectedness   + complex homogeneity under a compact group}
\\[-1pt]
&\hspace{35mm}\Longrightarrow H^0(M,\Omega_M^1)=0,\\[3pt]
&\nc\Tc=0,\quad \Tc\ne0,\quad M\text{ simply connected}
\\[-1pt]
&\hspace{35mm}\Longrightarrow H^0(M,\Omega_M^1)\ne0.
\end{aligned}
\end{equation}

\subsection{Holomorphic one-forms under a compact transitive group}
\label{sec:one-forms}

The first result asserts that any holomorphic $1$-form must be $d$-closed. This implies the vanishing 
of the holomorphic $1$-form under the assumption of simply-connectedness or the 
 finiteness of the fundamental group. 

\begin{lemma}\label{lem:one-forms}
Let $M$ be a compact connected complex manifold admitting a transitive
holomorphic action of a compact connected Lie group $K$. Every
holomorphic one-form on $M$ is closed. If $M$ is simply connected or with finite fundamental group, then
\begin{equation}\label{eq:h10-zero}
                         H^0(M,\Omega_M^1)=0.
\end{equation}
\end{lemma}

\begin{proof}
Put $\fk=\Lie(K)$. For $X\in\fk$, denote its real fundamental vector
field by $X^\#$. Since the action is holomorphic, $(X^\#)^{1,0}$ is a
holomorphic vector field. If $\alpha\in H^0(M,\Omega_M^1)$, then
\[
      \alpha(X^\#)=\alpha\bigl((X^\#)^{1,0}\bigr)
\]
is a holomorphic function on compact $M$, and hence is constant.
Complex-linear extension therefore defines
\[
 \ell_\alpha\in\fk_\C^*,\qquad
 \ell_\alpha(X)=\alpha(X^\#),\qquad \fk_\C=\fk\otimes_{\mathbb R}\C.
\]
The map $\alpha\mapsto\ell_\alpha$ is injective, since the fundamental
fields span the tangent space at every point. Its image
\[
 W=\{\ell_\alpha:\alpha\in H^0(M,\Omega_M^1)\}\subseteq\fk_\C^*
\]
is invariant under the coadjoint representation of $K$, by equivariance
of fundamental vector fields.

Fix $p\in M$. The evaluation map
\[
 \rho_p:\fk\longrightarrow T_pM,\qquad X\longmapsto X^\#_p
\]
is surjective. Its complexified dual is therefore injective, and
$\ell_\alpha=\rho_p^*\alpha_p$. Complex conjugation here and below is
defined by the real forms $\fk$ and $T_pM$. We have
\[
 W\subseteq\rho_p^*(T_p^{*1,0}M),\qquad
 \overline W\subseteq\rho_p^*(T_p^{*0,1}M).
\]
Since the two types of covectors have zero intersection,
\begin{equation}\label{eq:w-conjugate}
                            W\cap\overline W=\{0\}.
\end{equation}

Decompose the compact Lie algebra into its center and compact simple
ideals:
\[
             \fk=\mathfrak{z}\oplus\fk_1\oplus\cdots\oplus\fk_r.
\]
Its complex coadjoint representation is
\begin{equation}\label{eq:coadjoint-decomposition}
 \fk_\C^*=\mathfrak{z}_\C^*\oplus
            \bigoplus_{a=1}^{r}(\fk_a)_\C^*.
\end{equation}
The central summand is trivial. Each nontrivial summand is irreducible,
occurs with multiplicity one, and is preserved by complex conjugation.
These nontrivial summands are pairwise nonisomorphic as representations
of the full group $K$, even if some simple ideals are isomorphic as
abstract Lie algebras. Complete reducibility therefore implies that an
invariant subspace having a nonzero component in one of these
summands contains that entire summand. Applied to $W$, this would
contradict \eqref{eq:w-conjugate}. Consequently,
\[
                        W\subseteq\mathfrak{z}_\C^*,
\]
where $\mathfrak{z}_\C^*$ is identified with the annihilator of the semisimple
summand. In particular,
\begin{equation}\label{eq:annihilate-commutator}
                  \ell_\alpha([\fk,\fk])=0.
\end{equation}

For fundamental vector fields, constancy of the evaluations gives
\begin{align*}
 d\alpha(X^\#,Y^\#)
 &=X^\#\bigl(\alpha(Y^\#)\bigr)
   -Y^\#\bigl(\alpha(X^\#)\bigr)-\alpha([X^\#,Y^\#])\\
 &=-\alpha([X^\#,Y^\#])=0.
\end{align*}
The last equality follows from \eqref{eq:annihilate-commutator}; the
conventional sign relating $[X^\#,Y^\#]$ to $[X,Y]^\#$ does not affect the last quantity being zero.
Since the fundamental fields span $TM$, we obtain $d\alpha=0$.

If $M$ is simply connected, $\alpha=df$ for a global complex-valued
smooth function $f$. The form $\alpha$ has type $(1,0)$, so
$\bar\partial f=0$. Compactness makes the holomorphic function $f$
constant. Hence $\alpha=0$.
\end{proof}

\begin{remark}\label{rem:semisimple}
If $K$ is semisimple, then $\mathfrak{z}=0$, so the proof gives $W=0$ directly,
without using simple connectivity in Lemma~\ref{lem:one-forms}.
For $C$-spaces with compact semisimple transitive group, the vanishing
$h^{1,0}=0$ appeared first in 
\cite[Section~4.1, equation~(4.6)]{Lauret2026}.
\end{remark}

\subsection{Proof of the theorem and the homogeneous case}
\label{sec:conclusion}

\begin{lemma}\label{lem:fixed-vector}
Let $(M,J,g)$ be a connected Hermitian manifold with $\nc\Tc=0$.
For $p\in M$, put
\[
 V=T_p^{1,0}M,\qquad
 F=V^{\Hol_p^0(\nc)},
\]
where $\Hol_p^0(\nc)$ denotes the restricted holonomy group. Then
\begin{equation}\label{eq:fixed-vector}
                         F=0\quad\Longrightarrow\quad\Tc=0.
\end{equation}
Equivalently, nonzero parallel Chern torsion forces a nonzero fixed
vector of the restricted holonomy.
\end{lemma}

\begin{proof}
Write $\nabla=\nc$, $T=T^c$, and $\Rm=\Rm^c$. We may pass to the universal
cover: this does not change the restricted holonomy or the condition
$\nabla T=0$. Holonomy-invariant subspaces at $p$ can then be extended
to global parallel subbundles by parallel transport. The result is contained in Corollary \ref{cor:irreducible}.
\end{proof}

\begin{proof}[Proof of Theorem~\ref{thm:main}]
The identity component of a compact transitive group is transitive on
a connected manifold, so Lemma~\ref{lem:one-forms} applies and gives
\[
                         H^0(M,\Omega_M^1)=0.
\]
Suppose $\nc\Tc=0$ but $\Tc\ne0$. Lemma~\ref{lem:fixed-vector}
provides a nonzero vector
\[
                 u\in(T_p^{1,0}M)^{\Hol_p^0(\nc)}.
\]
Simple connectivity implies $\Hol_p(\nc)=\Hol_p^0(\nc)$. Thus $u$
extends by parallel transport to a global nonzero parallel
$(1,0)$-vector field $U$.

Define its Hermitian dual by
\begin{equation}\label{eq:parallel-form}
                         \alpha(X)=g(X,\overline U).
\end{equation}
The form $\alpha$ has type $(1,0)$, and metric compatibility gives
$\nc\alpha=0$. It is nonzero because
\[
                       \alpha(U)=g(U,\overline U)>0.
\]
On the holomorphic cotangent bundle, the $(0,1)$-part of the Chern
connection is the Dolbeault operator. Consequently,
\[
                 \nc\alpha=0\quad\Longrightarrow\quad
                 \bar\partial\alpha=0.
\]
This produces a nonzero element of $H^0(M,\Omega_M^1)$, a
contradiction. Hence $\Tc=0$. Vanishing Chern torsion is equivalent to
$d\omega=0$, where $\omega(X,Y)=g(JX,Y)$, so $g$ is K\"ahler.
\end{proof}

\begin{proof}[Proof of Corollary~\ref{cor:homogeneous}]
If the specified transitive group $G$ is compact, apply
Theorem~\ref{thm:main} directly. More generally, because $g$ is
$G$-invariant and the action is holomorphic, the image of $G$ lies in
the group of holomorphic isometries of $(M,g,J)$. This is a closed
subgroup of the compact isometry group of the compact Riemannian
manifold $(M,g)$. The closure $K$ of the image of $G$ is therefore a
compact transitive group of biholomorphisms. Theorem~\ref{thm:main}
again applies.
\end{proof}

\begin{remark} Once the underlying complex manifold admits a
compact transitive group of biholomorphisms, invariance of the metric
is unnecessary.
\end{remark}

\subsection*{Acknowledgments}
The author thanks J. Lauret and F. Pediconi for the discussions and suggestions. J. Lauret informed the author that 
Corollary \ref{cor:homogeneous} was also obtained by his student. ChatGPT was used in this project, 
particularly for some of the more elegant arguments.

\end{document}